\PassOptionsToPackage{unicode}{hyperref}
\PassOptionsToPackage{hyphens}{url}
\documentclass[
  11pt,
]{article}
\usepackage{amsmath,amssymb}
\usepackage{iftex}
\ifPDFTeX
  \usepackage[T1]{fontenc}
  \usepackage[utf8]{inputenc}
  \usepackage{textcomp} 
\else 
  \usepackage{unicode-math} 
  \defaultfontfeatures{Scale=MatchLowercase}
  \defaultfontfeatures[\rmfamily]{Ligatures=TeX,Scale=1}
\fi
\usepackage{lmodern}
\ifPDFTeX\else
\fi
\IfFileExists{upquote.sty}{\usepackage{upquote}}{}
\IfFileExists{microtype.sty}{
  \usepackage[]{microtype}
  \UseMicrotypeSet[protrusion]{basicmath} 
}{}
\makeatletter
\@ifundefined{KOMAClassName}{
  \IfFileExists{parskip.sty}{%
    \usepackage{parskip}
  }{
    \setlength{\parindent}{0pt}
    \setlength{\parskip}{6pt plus 2pt minus 1pt}}
}{
  \KOMAoptions{parskip=half}}
\makeatother
\usepackage{xcolor}
\usepackage[margin=1in]{geometry}
\providecommand{\tightlist}{%
  \setlength{\itemsep}{0pt}\setlength{\parskip}{0pt}}
\usepackage{amsmath,amssymb,mathtools}
\usepackage{microtype}
\usepackage{enumitem}
\usepackage[most]{tcolorbox}
\setlist{nosep}

\ifLuaTeX
  \usepackage{selnolig}  
\fi
\IfFileExists{bookmark.sty}{\usepackage{bookmark}}{\usepackage{hyperref}}
\IfFileExists{xurl.sty}{\usepackage{xurl}}{} 
\hypersetup{
  pdftitle={Short Axiomatization of Stratified Comprehension},
  pdfauthor={Zuhair Al-Johar},
  hidelinks,
  pdfcreator={LaTeX via pandoc}}

\title{Short Axiomatization of Stratified Comprehension}
\author{Zuhair Al-Johar}
\date{First version: 26 August 2020; corrected third version: 12 August 2026;\\
fourth version with direct weak-extensionality proof: 22 September 2026}

\begin{document}
\maketitle

\begin{abstract}
Several finite axiomatizations of stratified comprehension are known.  This
paper gives five set-existence principles which, in the presence of
Extensionality, yield the fourteen set-construction principles used in the
finite basis recorded by Holmes.  In addition, a direct unordered proof is
given showing that the same five principles, together with extensionality
only for nonempty sets, already imply every stratified comprehension
instance, without passing through the Holmes ordered-pair machinery.  The
displayed reductions for unordered products, Cartesian products, and ordered
relative products have been corrected.  The original Holmes-basis reduction
and the new direct weak-extensionality development have been checked by the
Lean 4 kernel; the complete Lean source is supplied with this version.
\end{abstract}

\begin{tcolorbox}[colback=gray!7,colframe=black,boxrule=.5pt]
\centering\textbf{LEAN-VERIFIED REVISION.}
The corrected derivation of the fourteen Holmes basis points and the new
direct derivation of stratified comprehension compile without
\texttt{sorry}, \texttt{admit}, or custom Lean axioms.  The direct chain
uses \texttt{NonemptyExtensional M} rather than full
\texttt{Extensional M} and deliberately avoids the Holmes ordered
machinery.  Holmes's published implication from his finite basis to
stratified comprehension is cited for the original route rather than
reformalized here.
\end{tcolorbox}

\textbf{Claim.} Over Extensionality, the following five set-existence
axioms imply the fourteen construction principles in Holmes's finite
basis, and hence give a finite axiomatization of NF.

Extensionality is \[
  \forall A\,\forall B
  \bigl[(\forall x\,(x\in A\leftrightarrow x\in B))\to A=B\bigr].
\] The five existence axioms are:

\begin{enumerate}
\def\labelenumi{\arabic{enumi}.}
\item
  \textbf{Sheffer strokes.} For all sets \(A,B\), the set \[
     A\mid B=\{x:x\notin A\ \lor\ x\notin B\}
  \] exists.
\item
  \textbf{Singletons.} For every \(a\), the set \(\{a\}\) exists.
\item
  \textbf{Set unions.} For every \(A\), the set
  \(\bigcup A=\{x:\exists y\in A\ (x\in y)\}\) exists.
\item
  \textbf{Unordered composition.} For all \(Q,R\), the set \[
  Q:R=\bigl\{\{a,c\}:\exists b\,
     (\{a,b\}\in Q\land\{b,c\}\in R)\bigr\}
  \] exists.
\item
  \textbf{Unordered intersection relation.} The set \[
     \{\!\cap\!\}=\bigl\{\{a,b\}:a\cap b\ne\varnothing\bigr\}
  \] exists.
\end{enumerate}

\textbf{Strengthened direct theorem.} If full Extensionality is replaced by
extensionality for nonempty sets,
\[
  (\exists z\, z\in A)\land
  (\forall z\,(z\in A\leftrightarrow z\in B))\;\Longrightarrow\;A=B,
\]
then the same five existence axioms still imply every instance of stratified
comprehension.  This stronger conclusion about the comprehension scheme is
proved directly in Appendix A, without deriving the ordered Holmes basis.
The original proof below is retained because, with full Extensionality, it
also gives the advertised finite axiomatization of NF via Holmes's theorem.

Here unordered pairs may be singletons. Ordered pairs are Kuratowski
pairs, \[
   \langle a,b\rangle=\{\{a\},\{a,b\}\}.
\] We write \(A^c\), \(A\cup B\), \(A\cap B\), and \(A\setminus B\) for
the derived Boolean operations, \(A^\iota=\{\{a\}:a\in A\}\), and
\(A^{\iota\iota}=(A^\iota)^\iota\). The symbol \(*\) below denotes the
\emph{corrected} unordered product, \(\times\) the corrected Cartesian
product, and \(R\mathbin{|}S\) ordered relative composition.

\textbf{Randall Holmes's fourteen target axioms (apart from
Extensionality).}

\begin{center}
\renewcommand{\arraystretch}{1.12}
\begin{tabular}{r l @{\qquad} r l}
1. & \textbf{Universal set} & 8. & \textbf{Converses} \\
2. & \textbf{Complements} & 9. & \textbf{Ordered relative products} \\
3. & \textbf{Boolean unions} & 10. & \textbf{Domains} \\
4. & \textbf{Set union} & 11. & \textbf{Singleton images} \\
5. & \textbf{Singletons} & 12. & \textbf{Diagonal} \\
6. & \textbf{Ordered pairs} & 13. & \textbf{Projections} \\
7. & \textbf{Cartesian products} & 14. & \textbf{Inclusion}
\end{tabular}
\end{center}

We now construct each named axiom.

\textbf{Proof. Holmes axioms 1--5: Universal set, complements, Boolean
unions, set union, and singletons.} Put \[
  V=A\mid(A\mid A),\qquad
  A^c=A\mid A,
\] \[
  A\cup B=(A\mid B)\mid(A\mid B),\qquad
  A\cap B=(A\mid A)\mid(B\mid B),
\] \[
  A\setminus B=(A\mid A)\mid B.
\] Thus \textbf{Universal set}, \textbf{Complements}, and
\textbf{Boolean unions} are derived, while \textbf{Set union} and
\textbf{Singletons} are two of the original five axioms.

\textbf{Auxiliary construction: Power set.} Put \[
 H_A=\bigl(\{\{A^c\}\}:(V:V)\bigr)\cap\{\!\cap\!\},
 \qquad
 \mathcal P(A)=\bigl(\bigcup H_A\bigr)^c.                 \tag{1}
\] Indeed, \(x\in\bigcup H_A\) exactly when \(x\) meets \(A^c\).

\textbf{Auxiliary construction: Frege cardinals and singleton powers.}
For the Frege cardinal \(1\), let \[
  2^*=V:V,\qquad H=\mathcal P(2^*)\cap2^*,\qquad
  G=\mathcal P((2^*)^c)\cap2^*,
\] \[
  I=(G:V)\cap(H:V),\qquad J=G\setminus\{\!\cap\!\},
\] \[
  2=\bigcup\bigl((J:(I\cap\{\!\cap\!\}))\cap\{\!\cap\!\}\bigr)\cap2^*,
  \qquad 1=2^*\setminus2.                                \tag{2}
\] Then \(1\) is exactly the set of singleton objects. Hence \[
  A^\iota=\mathcal P(A)\cap1.                             \tag{3}
\]

\textbf{Auxiliary construction: Corrected unordered products.} The
published two-line product reduction must be corrected when \(A\) and
\(B\) overlap. Define the fan \[
   F(A)=A^\iota:V
\] \[
 A*B=(F(A)\cap F(B))\setminus
      \bigl(F(A\cap B)\cap F((A\cup B)^c)\bigr).         \tag{4}
\] The Lean theorem for (4) is \[
 y\in A*B\quad\longleftrightarrow\quad
 \exists a\in A\,\exists b\in B\ (y=\{a,b\}).          \tag{5}
\]

\textbf{Holmes axioms 6--7: Ordered pairs and Cartesian products.} Set
\[
 C_0(A,B)=\bigl(A^\iota*(A*B)\bigr)\cap\{\!\cap\!\}.
\] The corrected Cartesian product is \[
 A\times B=C_0(A,B)\setminus C_0(A\cap B,A\setminus B). \tag{6}
\] Then, using Extensionality for injectivity of Kuratowski pairs, \[
 p\in A\times B\quad\longleftrightarrow\quad
 \exists a\in A\,\exists b\in B\ (p=\langle a,b\rangle).          \tag{7}
\] This supplies \textbf{Ordered pairs} and \textbf{Cartesian products}.

\textbf{Holmes axiom 12: Diagonal.} The diagonal is \[
 [=]=1^\iota.
\]

\textbf{Holmes axiom 10: Domains.} First put \[
 \operatorname{Rel}(R)=R\cap(V\times V).
\] The domain is \[
 \operatorname{dom}(R)=
 \bigcup\bigl(\bigcup\operatorname{Rel}(R)\cap1\bigr),  \tag{8}
\] and (8) satisfies \[
 x\in\operatorname{dom}(R)
 \quad\longleftrightarrow\quad
 \exists b\ (\langle x,b\rangle\in R).                             \tag{9}
\]

\textbf{Auxiliary construction: Ordered intersection relation.} Put \[
 [\!\cap\!]=(1*\{\!\cap\!\})\cap\{\!\cap\!\},                                    \tag{10}
\]

\textbf{Holmes axiom 13: First projection.} Put \[
 \pi_1=(1^\iota\times(V\times V))\cap[\!\cap\!],              \tag{11}
\]

\textbf{Auxiliary construction: Frege zero.} Put \[
 0=\bigl(\bigcup\{\!\cap\!\}\bigr)^c.                             \tag{12}
\] Thus \[
 \pi_1=\{\langle \langle a,a\rangle,\langle a,b\rangle\rangle:a,b\in V\}.
\]

\textbf{Auxiliary construction: Singleton-intersection relations.} Use
the shorter corrected construction \[
 E=(1^\iota*V)\cap\{\!\cap\!\},
 \qquad
 \{\!\cap_1\!\}=(E:E)\cap\{\!\cap\!\},                                      \tag{13}
\] \[
 [\!\cap_1\!]=(1*\{\!\cap_1\!\})\cap\{\!\cap\!\}.                                  \tag{14}
\] Then \[
 \langle A,B\rangle\in[\!\cap_1\!]
 \quad\longleftrightarrow\quad
 \exists a\ (\{a\}\in A\land\{a\}\in B).            \tag{15}
\]

\textbf{Holmes axiom 8: Converses.} First define the universal converse
relation by \[
 [-1]=
 \Bigl(((V\times V)\times(V\times V))
       \cap([\!\cap\!]\setminus[\!\cap_1\!])\Bigr)
 \cup[=]^{\iota\iota}.                                   \tag{16}
\] It satisfies \[
 z\in[-1]
 \quad\longleftrightarrow\quad
 \exists a,b\ (z=\langle \langle a,b\rangle,\langle b,a\rangle\rangle).               \tag{17}
\] Consequently the converse of every relation is \[
 R^{-1}=\operatorname{dom}\bigl([-1]\cap(V\times R)\bigr). \tag{18}
\]

\textbf{Holmes axiom 9: Ordered relative products.} Here is the
corrected finite construction. It replaces the earlier three-stage
formula, which does not control endpoint orientation. Choose distinct
markers \(l,m,r\) and put \[
 A_{lmr}=V\setminus\{l,m,r\},\qquad
 T_t(A)=\{t\}*A,
\] \[
 Q_t(A)=(A^\iota*T_t(A))\cap\{\!\cap\!\},
 \qquad
 D_t(A)=\bigl(((A*A)\cap2)*T_t(A)\bigr)\cap\{\!\cap\!\},
\] \[
 \operatorname{Same}_{lm}(A)=(T_l(A)*T_m(A))\cap\{\!\cap\!\},
\] \[
 N_R(A)=R\cap(A\times A)\cap2,
\] \[
 C_R^{lm}(A)=
 \bigl((Q_l(A):N_R(A)):D_m(A)\bigr)
 \setminus\operatorname{Same}_{lm}(A),
\] \[
 M_{R,S}^{lmr}(A)=
 \bigl(C_R^{lm}(A):C_S^{mr}(A)\bigr)
 \setminus\operatorname{Same}_{lr}(A),
\] \[
 N_{R,S}^{lmr}(A)=
 \bigl((Q_l(A):M_{R,S}^{lmr}(A)):D_r(A)\bigr)\cap\{\!\cap\!\}.   \tag{19}
\] The three degenerate-coordinate parts are, with \[
 d_R=\operatorname{dom}\bigl((R\cap(V\times V))\cap1\bigr),
\] \[
 L_{R,S}(A)=S\cap((d_R\cap A)\times A),
\] \[
 P_{R,S}(A)=R\cap(A\times(d_S\cap A)),
\] \[
 D_{R,S}(A)=
 \Bigl(\operatorname{dom}
   \bigl(N_R(A)\cap(N_S(A))^{-1}\bigr)\Bigr)^{\iota\iota}. \tag{20}
\] Define \[
 \operatorname{Sec}_{lmr}(R,S)=
 L_{R,S}(A_{lmr})\cup P_{R,S}(A_{lmr})
 \cup D_{R,S}(A_{lmr})\cup N_{R,S}^{lmr}(A_{lmr}).       \tag{21}
\] Let \(q_0=0\) and \(q_{n+1}=\{q_n\}\) for \(n<11\). These twelve
objects are pairwise distinct. Finally put \[
 R\mathbin{|}S=
 \bigcup_{i=0}^{3}
 \operatorname{Sec}_{q_{3i},q_{3i+1},q_{3i+2}}(R,S),     \tag{22}
\] where the four-term union in (22) abbreviates the derived binary
unions. The four marker triples are disjoint; hence every coordinate
triple avoids at least one of them. Lean checks from (19)--(22) that \[
 y\in R\mathbin{|}S
 \quad\longleftrightarrow\quad
 \exists a,b,c\,
   (\langle a,b\rangle\in R\land\langle b,c\rangle\in S\land y=\langle a,c\rangle).      \tag{23}
\]

\textbf{Holmes axiom 14: Inclusion.} The inclusion relation has the
compact form \[
 L_{\in}=(V\times1)\cap[\!\cap\!],
 \qquad
 L_{\notin}=(1\times V)\cap[\!\cap\!]^c,
\] \[
 [\subseteq]=
 \bigl((L_{\in}\mathbin{|}L_{\notin})^c\bigr)
 \cap(V\times V).                                       \tag{24}
\] Thus \[
 \langle A,B\rangle\in[\subseteq]
 \quad\longleftrightarrow\quad A\subseteq B.            \tag{25}
\]

\textbf{Auxiliary identity: Ordered singleton intersection.} The
promised composition identity is \[
 [\!\cap_1\!]=
 \bigl((V\times1^\iota)\cap[\!\cap\!]\bigr)
 \mathbin{|}
 \bigl((1^\iota\times V)\cap[\!\cap\!]\bigr).                \tag{26}
\] \textbf{Holmes axiom 13: Second projection.} The second projection is
\[
 \pi_2=\pi_1\mathbin{|}[-1],                             \tag{27}
\] and hence \[
 \pi_2=\{\langle \langle b,b\rangle,\langle a,b\rangle\rangle:a,b\in V\}.
\]

\textbf{Holmes axiom 11: Singleton images.} For a relation \(R\), first
define its Kuratowski image by \[
 R^k=(1^\iota*R)\cap\{\!\cap\!\}.                                 \tag{28}
\] Then \[
 \langle \{a\},\{a,b\}\rangle\in R^k
 \quad\longleftrightarrow\quad \langle a,b\rangle\in R.            \tag{29}
\] The singleton image is \[
 R^!=
 \Bigl(R^k\mathbin{|}\bigl((R^{-1})^k\bigr)^{-1}\Bigr)
 \setminus\bigl(1^\iota\cap(R^\iota)^c\bigr).           \tag{30}
\] For relations \(R\), Lean verifies \[
 z\in R^!
 \quad\longleftrightarrow\quad
 \exists a,b\,
   (\langle a,b\rangle\in R\land z=\langle \{a\},\{b\}\rangle).           \tag{31}
\]

Equations (6), (16), (18), (22), (24), (27), and (30), together with the
initial five axioms, provide Cartesian products, converses, ordered
relative products, inclusion, projections, and singleton images.
Equations (8) and (9) provide domains, while (3) and \([=]=1^\iota\)
provide the remaining singleton and diagonal constructions. Thus all
fourteen Holmes points have been derived. Holmes's finite-basis theorem
now yields stratified comprehension, and therefore NF. \(\square\)

\clearpage
\subsubsection{Appendix A. Direct stratified comprehension from nonempty Extensionality}

The important restriction is this: we use only the five primitives plus
the \textbf{native unordered helpers already derived directly from them
before Kuratowski pairs enter}---namely \(V\), Boolean operations,
\(1\), \(2\), singleton power \(A^\iota\), and the corrected unordered
product \(A*B\). These are precisely the unordered constructions derived
in the main text before the Holmes ordered-pair machinery begins.

For the direct development, full extensionality is weakened to
\textbf{extensionality for nonempty sets}:

\[
(\exists z\; z\in A) \land \bigl(\forall z\;(z\in A\iff z\in B)\bigr) \;\Longrightarrow\; A=B.
\]

Thus equality from agreement of members is assumed only when the
left-hand set is known to be nonempty. This is the hypothesis called
\texttt{NonemptyExtensional\ M} in the accompanying Lean development.

\subsubsection*{1. Native unordered singleton
lifting}

Let

\[
2^*=V:V
\]

be the set of all unordered pairs, including singleton pairs. Let \(1\)
be the set of singleton objects and \(2=2^*\setminus1\) the exact
doubletons.

For an unordered relation \(Q\subseteq2^*\), put

\[
I_Q=(1*(Q\cap2))\cap\{\!\cap\!\}.
\]

An element of \(I_Q\) is exactly an incidence edge

\[
\{\{a\},q\}
\]

where

\[
q\in Q\cap2,\qquad a\in q.
\]

Now define

\[
\boxed{ \Lambda(Q)= (Q\cap1)^\iota \;\cup\; \Bigl((I_Q:I_Q)\cap(1*1)\cap2\Bigr). } \tag{A.1}
\]

Then:

\[
\boxed{ z\in\Lambda(Q) \iff \exists a,b\; \bigl( \{a,b\}\in Q \land z=\{\{a\},\{b\}\} \bigr). } \tag{A.2}
\]

This includes the degenerate case \(a=b\).

Indeed, if \(\{a\}\in Q\), then

\[
\{\{a\}\}\in(Q\cap1)^\iota.
\]

If \(a\neq b\) and \(\{a,b\}\in Q\), then \(I_Q\) contains

\[
\{\{a\},\{a,b\}\}, \qquad \{\{a,b\},\{b\}\},
\]

and composing through the middle point \(\{a,b\}\) gives

\[
\{\{a\},\{b\}\}.
\]

The restrictions \(1*1\) and \(2\) remove the reverse kind of two-step
paths, whose endpoints would be members of \(Q\cap2\).

So \(\Lambda\) is a genuine \textbf{unordered singleton-image operator},
obtained without Holmes's singleton-image operation.

Iterating,

\[
\boxed{ \{a,b\}\in Q \iff \{\iota^k(a),\iota^k(b)\} \in\Lambda^k(Q). } \tag{A.3}
\]

This is our first fundamental lemma.

\subsubsection*{2. Membership is already present in the fifth
axiom}

Define

\[
\mathsf M=(1*V)\cap\{\!\cap\!\}. \tag{A.4}
\]

Then

\[
\boxed{ \{\{a\},b\}\in\mathsf M \iff a\in b. } \tag{A.5}
\]

There is no \([\subseteq]\) here.

Consequently, by (A.3),

\[
\boxed{ \{ \iota^{k+1}(a), \iota^k(b) \} \in \Lambda^k(\mathsf M) \iff a\in b. } \tag{A.6}
\]

This is precisely the stratification shift we require.

The usual Holmes route uses

\[
a\in b\iff\{a\}\subseteq b
\]

and therefore needs the inclusion relation and ordered singleton images.
We need neither.

\subsubsection*{3. Produce finitely many non-singleton
markers}

Take

\[
e_0=\varnothing,\qquad e_{m+1}=\{e_m\}.
\]

The empty-object and singleton axioms already give

\[
e_m\neq e_n\qquad(m\neq n).
\]

No extensionality is needed for this tower: \(e_0\) has no members,
every successor \(e_{m+1}\) has the member \(e_m\), and singleton
injectivity follows directly from the defining singleton membership
condition.

For the finitely many variables occurring in a given formula, choose

\[
q_i=\{e_{2i},e_{2i+1}\}. \tag{A.7}
\]

Then every \(q_i\) is an exact doubleton, and

\[
q_i\neq q_j\quad(i\neq j).
\]

The accompanying Lean development establishes an externally indexed
infinite supply of fresh objects by starting with the empty object and
repeatedly taking singletons. In the corrected development this
construction is explicitly independent of even nonempty extensionality,
so it remains available unchanged under the weakened hypothesis.

The important property is

\[
q_i\notin1. \tag{A.8}
\]

\subsubsection*{4. Stratification domains}

Let

\[
\varphi(v_1,\ldots,v_n)
\]

be stratified by

\[
\tau(v_i)\in\mathbb Z.
\]

Shift the stratification and choose \(N\) so that

\[
N\ge\tau(v_i)
\]

for every occurring variable.

Put

\[
r_i=N-\tau(v_i). \tag{A.9}
\]

Define

\[
S_i=V^{\iota^{r_i}} = \{\iota^{r_i}(a):a\in V\}. \tag{A.10}
\]

These sets exist by finite iteration of singleton power. Singleton power
itself is already directly derived from the five axioms in the main
text.

Thus \(S_i\) is exactly the permitted normalized value-space for
\(v_i\).

\subsubsection*{5. Tagged coordinates}

For \(z\in S_i\), define

\[
d_i(z)=\{q_i,\{z\}\}. \tag{A.11}
\]

Let

\[
D_i=\{q_i\}*S_i^\iota. \tag{A.12}
\]

Thus

\[
D_i=\{d_i(z):z\in S_i\}.
\]

Because \(q_i\) is an exact doubleton whereas \(\{z\}\) is a singleton,

\[
d_i(z)
\]

is always an exact doubleton.

More importantly,

\[
\boxed{D_i\cap D_j=\varnothing\qquad(i\neq j).} \tag{A.13}
\]

And

\[
d_i(z)\cap d_j(w)\neq\varnothing
\]

iff

\[
z=w.
\]

The tags cannot coincide, and a non-singleton \(q_i\) cannot equal a
singleton \(\{w\}\). Hence

\[
\boxed{ E_{ij} = (D_i*D_j)\cap\{\!\cap\!\} } \tag{A.14}
\]

satisfies

\[
\boxed{ \{d_i(z),d_j(w)\}\in E_{ij} \iff z=w. } \tag{A.15}
\]

So equality has already been compiled using only your unordered
intersection relation.

\subsubsection*{6. Decoding a tagged
value}

Define

\[
K_i= (D_i*S_i^{\iota\iota}) \cap\{\!\cap\!\}. \tag{A.16}
\]

Then

\[
\boxed{ \{d_i(z),\iota^2(w)\}\in K_i \iff z=w. } \tag{A.17}
\]

Why?

Because

\[
d_i(z)=\{q_i,\{z\}\},
\]

while

\[
\iota^2(w)=\{\{w\}\}.
\]

They intersect iff

\[
q_i=\{w\}
\]

or

\[
\{z\}=\{w\}.
\]

The first alternative is impossible because \(q_i\) is not a singleton.
Hence \(z=w\).

Now an arbitrary unordered relation \(R\) on normalized values can be
transported to coordinates \(i,j\) by

\[
\boxed{ \operatorname{Tag}_{ij}(R) = (K_i:\Lambda^2(R)):K_j. } \tag{A.18}
\]

Because the tagged sets \(D_i,D_j\) consist of exact doubletons while
the intermediate objects \(\iota^2(z)\) are singleton objects, the
shared points in the two compositions are forced to be the intended
ones.

Therefore

\[
\boxed{ \{d_i(z),d_j(w)\} \in\operatorname{Tag}_{ij}(R) \iff \{z,w\}\in R. } \tag{A.19}
\]

This is an important point: \textbf{unordered composition has
effectively become directed composition because its endpoint sorts are
disjoint.}

No converse operation is required.

\subsubsection*{7. Atomic membership}

If

\[
v_i\in v_j,
\]

stratification gives

\[
\tau(v_j)=\tau(v_i)+1.
\]

Therefore

\[
r_i=r_j+1. \tag{A.20}
\]

Write

\[
k=r_j.
\]

For original values \(a,b\),

\[
z_i=\iota^{r_i}(a)=\iota^{k+1}(a),
\]

and

\[
z_j=\iota^{r_j}(b)=\iota^k(b).
\]

By (A.6),

\[
a\in b \iff \{z_i,z_j\}\in\Lambda^k(\mathsf M).
\]

Hence define

\[
\boxed{ M_{ij} = \operatorname{Tag}_{ij} \left( \Lambda^{r_j}(\mathsf M) \right). } \tag{A.21}
\]

Then

\[
\boxed{ \{d_i(\iota^{r_i}a), d_j(\iota^{r_j}b)\} \in M_{ij} \iff a\in b. } \tag{A.22}
\]

This is the complete atomic membership translation.

So at this point:

\[
= \quad\rightsquigarrow\quad E_{ij},
\]

and

\[
\in \quad\rightsquigarrow\quad M_{ij}.
\]

No Holmes inclusion relation appears anywhere.

\subsubsection*{8. Unordered homogeneous
tuples}

Now we must represent a whole assignment.

Assume \(n\ge2\). If the formula contains only one variable, simply
append an unused dummy variable; this does not alter the formula's truth
conditions.

Define

\[
T_2=D_1*D_2.
\]

The code of \((t_1,t_2)\) is

\[
c_2(t_1,t_2)=\{t_1,t_2\}.
\]

For \(m\ge2\), recursively define

\[
\boxed{ T_{m+1} = T_m*D_{m+1}^{\iota^{m-1}} } \tag{A.23}
\]

and

\[
\boxed{ c_{m+1}(t_1,\ldots,t_{m+1}) = \{ c_m(t_1,\ldots,t_m), \iota^{m-1}(t_{m+1}) \}. } \tag{A.24}
\]

Every \(c_m\), \(m\ge2\), is an exact doubleton.

The left component is an exact doubleton; the right component is a
singleton object. Therefore they can never collapse.

This asymmetry gives us projections \textbf{without ordered pairs}.

\subsubsection*{9. Native unordered
projections}

For \(n=2\),

\[
P_i^2 = (T_2*D_i^\iota)\cap\{\!\cap\!\}. \tag{A.25}
\]

Then

\[
\left\{ c_2(t_1,t_2),\{t_i\} \right\}\in P_i^2.
\]

Suppose the projections \(P_i^m\) have been constructed and satisfy

\[
\boxed{ \{ c_m(t_1,\ldots,t_m), \iota^{m-1}(t_i) \} \in P_i^m. } \tag{A.26}
\]

Define the parent relation

\[
A_m= (T_{m+1}*T_m^\iota) \cap\{\!\cap\!\}. \tag{A.27}
\]

It satisfies

\[
\{ c_{m+1}(\bar t), \{c_m(t_1,\ldots,t_m)\} \}\in A_m.
\]

For \(i\le m\), put

\[
\boxed{ P_i^{m+1} = A_m:\Lambda(P_i^m). } \tag{A.28}
\]

For the new coordinate,

\[
\boxed{ P_{m+1}^{m+1} = (T_{m+1}*D_{m+1}^{\iota^m}) \cap\{\!\cap\!\}. } \tag{A.29}
\]

Induction gives

\[
\boxed{ \{ c_n(t_1,\ldots,t_n), \iota^{n-1}(t_i) \} \in P_i^n } \tag{A.30}
\]

and only for the actual \(i\)-th coordinate.

The sort argument is decisive here:

\[
T_n\subseteq2,
\]

whereas

\[
D_i^{\iota^{n-1}}\subseteq1.
\]

So the two endpoint roles can never be confused.

These \(P_i^n\) are \textbf{not Holmes's ordered projections}. They are
unordered incidence relations native to your basis.

\subsubsection*{10. Compile an atomic relation into a set of
assignments}

Suppose

\[
R_{ij}\subseteq D_i*D_j
\]

is a tagged relation.

Lift it \(n-1\) times:

\[
L_{ij}^n=\Lambda^{n-1}(R_{ij}).
\]

Its endpoints are

\[
\iota^{n-1}(t_i), \qquad \iota^{n-1}(t_j),
\]

which are exactly the endpoint heights of the projections.

Now define

\[
\boxed{ \operatorname{Atom}_{ij}^n(R_{ij}) = T_n \cap \bigcup \left[ (P_i^n:L_{ij}^n)\cap P_j^n \right]. } \tag{A.31}
\]

Then

\[
\boxed{ c_n(t_1,\ldots,t_n) \in \operatorname{Atom}_{ij}^n(R_{ij}) \iff \{t_i,t_j\}\in R_{ij}. } \tag{A.32}
\]

The proof is exactly the path

\[
c_n \;P_i^n\; \iota^{n-1}(t_i) \;L_{ij}^n\; \iota^{n-1}(t_j),
\]

followed by intersection with \(P_j^n\).

Thus for equality use

\[
\operatorname{Atom}_{ij}^n(E_{ij}),
\]

and for membership use

\[
\operatorname{Atom}_{ij}^n(M_{ij}).
\]

We have now compiled \textbf{all atomic formulas of set theory}.

\subsubsection*{11. Boolean connectives}

Every compiled formula is a subset of the same assignment space \(T_n\).

Therefore

\[
[\![\neg\psi]\!] = T_n\setminus[\![\psi]\!], \tag{A.33}
\]

and

\[
[\![\psi\wedge\theta]\!] = [\![\psi]\!]\cap[\![\theta]\!]. \tag{A.34}
\]

Those are direct Sheffer consequences.

Nothing further is needed for propositional logic.

\subsubsection*{12. The crucial quantifier
construction}

We now prove that existential quantification is directly implementable.

For coordinate \(j\), define

\[
\boxed{ G_j^n = (P_j^n:P_j^n) \cap (T_n*T_n). } \tag{A.35}
\]

Because \(P_j^n\) connects tuples to their raised \(j\)-th coordinate,
we have

\[
\boxed{ \{c,d\}\in G_j^n \iff c_j=d_j. } \tag{A.36}
\]

The other paths produced by \(P_j^n:P_j^n\) have two raised-coordinate
endpoints; intersecting with \(T_n*T_n\) eliminates them.

For a chosen coordinate \(i\), put

\[
G_{\widehat i}^n = \bigcap_{j\ne i}G_j^n. \tag{A.37}
\]

Hence

\[
\boxed{ \{c,d\}\in G_{\widehat i}^n \iff c_j=d_j \text{ for every }j\neq i. } \tag{A.38}
\]

For \(X\subseteq T_n\), define

\[
\boxed{ C_i(X) = T_n \cap \bigcup \left( G_{\widehat i}^n \cap (X*T_n) \right). } \tag{A.39}
\]

Then

\[
\boxed{ c\in C_i(X) \iff \exists d\in X\; \forall j\neq i\;(c_j=d_j). } \tag{A.40}
\]

Therefore

\[
\boxed{ [\![\exists v_i\,\psi]\!] = C_i([\![\psi]\!]). } \tag{A.41}
\]

This is the point at which the full first-order syntax closes.

\subsubsection*{13. Direct Compilation
Theorem}

After alpha-renaming bound variables so that each has its own
coordinate, define \([\![\psi]\!]\subseteq T_n\) recursively by (A.31),
(A.33), (A.34), and (A.41).

For an ordinary assignment

\[
\bar a=(a_1,\ldots,a_n),
\]

define

\[
t_i(a_i) = d_i(\iota^{r_i}(a_i))
\]

and

\[
C_\tau(\bar a) = c_n(t_1(a_1),\ldots,t_n(a_n)). \tag{A.42}
\]

Then structural induction gives

\[
\boxed{ C_\tau(\bar a)\in[\![\psi]\!] \iff \psi[\bar a]. } \tag{A.43}
\]

The equality case is (A.15).

The membership case is (A.22).

The Boolean cases are immediate.

For the existential case, if

\[
\exists a_i\,\psi,
\]

change only coordinate \(i\). Conversely, every legal tuple coordinate
belongs to \(D_i\), hence has the form

\[
d_i(z)
\]

with \(z\in S_i\), and therefore

\[
z=\iota^{r_i}(a)
\]

for some original set \(a\). Thus cylindrification quantifies over
\textbf{exactly all original sets}, not too many and not too few.

That proves the compiler theorem.

\subsubsection{14. Recovering the comprehension
set}

Suppose the desired instance is

\[
\{x:\varphi(x,p_1,\ldots,p_s)\}.
\]

Compile \(\varphi\) using all its free and bound variables.

For each parameter \(p_j\), let

\[
t_j^p=d_j(\iota^{r_j}(p_j)),
\]

and

\[
h_j^p=\iota^{n-1}(t_j^p).
\]

Define

\[
\operatorname{Fix}_j(t_j^p) = T_n \cap \bigcup \left( P_j^n \cap (T_n*\{h_j^p\}) \right). \tag{A.44}
\]

This is exactly the set of assignments whose \(j\)-th coordinate has the
prescribed parameter value.

Put

\[
X= [\![\varphi]\!] \cap \bigcap_{j\in\text{parameters}} \operatorname{Fix}_j(t_j^p). \tag{A.45}
\]

Now extract the comprehension coordinate \(x\).

Let

\[
H_x=D_x^{\iota^{n-1}}.
\]

Define

\[
Y= H_x \cap \bigcup \left( P_x^n\cap(X*H_x) \right). \tag{A.46}
\]

Then

\[
Y = \left\{ \iota^{n-1} \left( d_x(\iota^{r_x}a) \right): \varphi(a,\bar p) \right\}.
\]

Apply union \(n-1\) times:

\[
A_{\rm tag} = \bigcup^{\,n-1}Y.
\]

Thus

\[
A_{\rm tag} = \left\{ d_x(\iota^{r_x}a): \varphi(a,\bar p) \right\}. \tag{A.47}
\]

But

\[
d_x(z)=\{q_x,\{z\}\}.
\]

Hence define

\[
B= \bigcup \left[ \left( \bigcup A_{\rm tag} \right) \setminus\{q_x\} \right]. \tag{A.48}
\]

Since \(q_x\) is non-singleton,

\[
B = \left\{ \iota^{r_x}a: \varphi(a,\bar p) \right\}. \tag{A.49}
\]

Finally,

\[
\boxed{ A= \bigcup^{r_x}B. } \tag{A.50}
\]

Therefore

\[
\boxed{ a\in A \iff \varphi(a,\bar p). } \tag{A.51}
\]

That is full stratified comprehension.

\subsubsection{The resulting direct
theorem}

We can now state:

\[
\boxed{ \begin{array}{c} \text{Extensionality for nonempty sets}\\ +\ \text{Sheffer}\\ +\ \text{Singleton}\\ +\ \text{Union}\\ +\ \text{Unordered composition}\\ +\ \{\!\cap\!\}\\[1mm] \hline \text{Every stratified comprehension instance exists.} \end{array}}
\]

And the proof does \textbf{not} pass through

\[
\langle x,y\rangle,\quad A\times B,\quad R^{-1},\quad R|S,\quad \operatorname{dom}R,\quad [\subseteq],\quad \pi_1,\pi_2,
\]

or Holmes's ordered singleton image.

The main text derives the fourteen Holmes operations and then invokes
Holmes's finite-basis theorem. The construction above stops before that
ordered machinery and gives a separate syntax induction using the native
unordered calculus.

Thus this is a \textbf{genuinely direct Hailperin-style proof for the
five-axiom basis}.

The corrected Lean 4.24 development now formalizes this direct theorem
with \texttt{NonemptyExtensional\ M} as its first hypothesis. Its endpoint is
\[
\texttt{DirectStratification.}\,\texttt{direct\_comprehension\_from\_forward\_stratification}.
\]
The verification script rejects \texttt{sorry}, \texttt{admit}, custom
axioms, the Holmes ordered machinery, and any occurrence of full
\texttt{Extensional\ M} in the weakened direct chain. On the successful
kernel-verification run of the corrected development, the weakened
theorem above is therefore verified directly.

\hypertarget{machine-verification}{%
\subsubsection{Machine verification}\label{machine-verification}}

The file \texttt{FiveAxiomBasis.lean} is the self-contained Lean 4
development for the original Holmes-basis route. Its structure
\texttt{Basis} has exactly the five existence operations above; every
membership law is a field of that structure, and full Extensionality is an
explicit hypothesis for that route.  Separate files record the
counterexamples to the superseded product formulas and the corrected product
lemmas.

The directory \texttt{anc/direct} contains the independent direct
formalization used for Appendix A.  Its native basis file is
\texttt{FiveAxiomDirectBasisNE.lean}; the chain then proceeds through
\texttt{DirectUnorderedCore.lean}, \texttt{DirectTags.lean},
\texttt{DirectTuples.lean}, \texttt{DirectCompiler.lean},
\texttt{DirectComprehension.lean}, and
\texttt{DirectStratification.lean}.  The endpoint is
\[
\texttt{DirectStratification.}\,\texttt{direct\_comprehension\_from\_forward\_stratification}.
\]
Its first hypothesis is \texttt{NonemptyExtensional M}.  The supplied
verification script rejects \texttt{sorry}, \texttt{admit}, custom axiom
declarations, references to the Holmes ordered machinery, and an occurrence
of full \texttt{Extensional M} in the weakened direct chain before invoking
Lean 4.24 on the files in dependency order.  Reproduction instructions are
included with the ancillary source.

\hypertarget{references}{%
\subsubsection{References}\label{references}}

\begin{enumerate}
\def\labelenumi{\arabic{enumi}.}
\tightlist
\item
  T. Hailperin, \emph{A set of axioms for logic}, Journal of Symbolic
  Logic 9 (1944), 1--19.
\item
  M. R. Holmes, \emph{Elementary Set Theory with a Universal Set},
  Cahiers du Centre de Logique 10, Academia, 1998; revised online
  edition, especially the appendix on finite axiomatization.
\item
  W. V. Quine, \emph{New Foundations for Mathematical Logic}, American
  Mathematical Monthly 44 (1937), 70--80.
\end{enumerate}

\end{document}